\documentclass[11pt]{amsart}
\usepackage{float}
\usepackage[usenames]{color}
\usepackage{graphicx}
\usepackage{amssymb}
\usepackage{amscd}
\usepackage{makecell}
\theoremstyle{plain}
\newtheorem{thm}{Theorem}[section]
\newtheorem{lemma}[thm]{Lemma}

\theoremstyle{definition}

\newcommand{\fra}{\mathfrak{a}}

\newcommand{\frg}{\mathfrak{g}}
\newcommand{\frh}{\mathfrak{h}}

\newcommand{\frk}{\mathfrak{k}}

\newcommand{\frp}{\mathfrak{p}}

\newcommand{\frt}{\mathfrak{t}}

\newcommand{\bbH}{\mathbb{H}}

\newcommand{\bbR}{\mathbb{R}}

\begin{document}

\title[Non-decreasable $K$-types are unitarily small]
{Non-decreasable $K$-types are unitarily small}

\author{Chao-Ping Dong}
\address[Dong]{School of Mathematical Sciences, Soochow University, Suzhou 215006,
P.~R.~China}
\email{chaopindong@163.com}

\author{Chengyu Du}
\address[Du]{School of Mathematical Sciences, Soochow University, Suzhou 215006,
P.~R.~China}
\email{cydu0973@suda.edu.cn}

\author{Haojun Xu}
\address[Xu]{School of Mathematical Sciences, Soochow University, Suzhou 215006,
P.~R.~China}
\email{20234207003@stu.suda.edu.cn}

\abstract{Let $G$ be a connected simple non-compact real reductive Lie group with a maximal compact subgroup $K$. This note aims to show that any non-decreasable $K$-type (in the sense of the first named author) is unitarily small (in the sense of Salamanca-Riba and Vogan). This answers Conjecture 2.1 of \cite{D} in the affirmative.}
 \endabstract



\maketitle
\section{Introduction}
Let $G$ be a connected simple non-compact Lie group in the Harish-Chandra class (see Section 3 of \cite{HC}).
Choose a Cartan involution $\theta$ such that $K:=G^{\theta}$ is a maximal compact subgroup of $G$. Let $\frg_0$ (resp. $\frk_0$) be the Lie algebra of $G$ (resp. $K$). Write
$$
\frg_0=\frk_0 + \frp_0
$$
as the Cartan decomposition on the Lie algebra level. Let $T$ be a maximal torus  of $K$, and we identify the abelian group $\widehat{T}$ of characters of $T$ with a lattice in $i\frt_0^*$. Here $\frt_0$ stands for the real Lie algebra of $T$. Let $\fra_0$ be the centralizer of $\frt_0$ in $\frp_0$. Then
$$
\frh_0:=\frt_0+\fra_0
$$ is a $\theta$-stable Cartan subalgebra of $\frg_0$. By a \textit{$K$-type}, we mean an irreducible representation of $K$. Note that any $K$-type must be finite-dimensional.  Fix a positive root system $\Delta^+(\frk, \frt)$ and denote the half sum of its members as $\rho_{\rm c}$.
Here we drop the subscript ``0" to denote the complexification of a Lie algebra.
We may and we will refer to a $K$-type by its highest weight. An important way to understand a representation $\pi$ of $G$ is to understand its $K$-types.

To reduce the workload of the classification of the \textit{Dirac series} of $G$, that is, all the irreducible unitary representations of $G$ having non-zero \textit{Dirac cohomology} (see \cite{Vog97} and \cite{HP} for the relevant backgrounds), the first named author proposed an approach to sharpen the Helgason-Johnson bound  in 1969 \cite{HJ}  and obtained explicit results for all exceptional Lie groups in \cite{D}. The improved Helgason-Johnson bound turned out to be quite helpful for  classifying the Dirac series of high rank exceptional Lie groups such as $E_{7(7)}$ \cite{DDW} and $E_{8(-24)}$ \cite{DDDLY}.

Choose any positive root system $\Delta^+(\frg, \frt)$ containing the fixed $\Delta^+(\frk, \frt)$. Write $\rho$ as the half sum of the roots in $\Delta^+(\frg, \frt)$, and put  $\rho_{\rm n}=\rho-\rho_{\rm c}$. The \textit{unitarily small} (\textit{u-small} for short henceforth) \textit{convex hull}, as defined by Salamanca-Riba and Vogan \cite{SV},  is the convex hull generated by the points $\{2w \rho_{\rm n} \mid w \in W(\mathfrak{g}, \mathfrak{t})\}$. A $K$-type is \textit{u-small} if its highest weight lives in the u-small convex hull. Otherwise, we will say that the $K$-type is \textit{u-large}.

\begin{figure}[H]
\centering
\scalebox{0.8}{\includegraphics{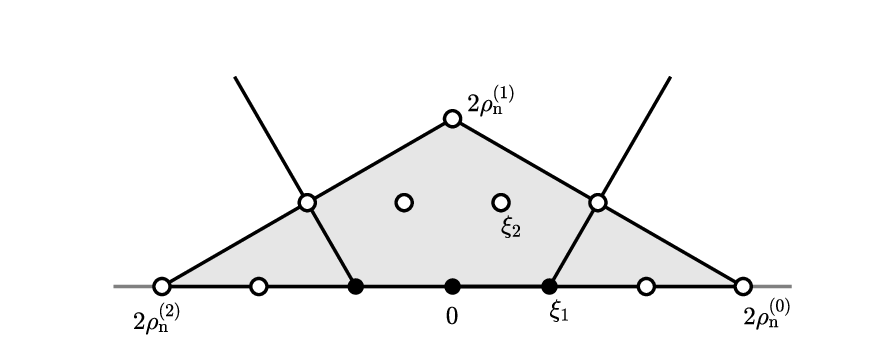}}
\caption{Non-decreasable $K$-types and u-small $K$-types for $U(2, 1)$}
\label{u21}
\end{figure}

In sharpening the Helgason-Johnson bound, u-small  $K$-types play a key role. One also needs to move a u-large $K$-type downwards, and this leads to the notion of non-decreasable $K$-types in Section 2.3 of \cite{D}.
Let $\Pi:=\{\alpha_i\mid 1\leqslant i\leqslant n\}$ be the simple roots of $\Delta^+(\frg, \frt)$, and let $\{\xi_i\}_{1\leqslant i\leqslant n}$ be the corresponding fundamental weights.
Let $\mu$ be the highest weight of a $K$-type such that $\mu+2\rho_{\rm c}$ is dominant for $\Delta^+(\frg, \frt)$.
Then $\mu$ is \textit{non-decreasable} if for any $1\leqslant i\leqslant n$, either the weight $\mu-\xi_i$ is not dominant for $\Delta^+(\mathfrak{k}, \mathfrak{t})$ or the weight $\mu+2\rho_{\rm c}-\xi_i$ is not dominant for $\Delta^+(\frg, \frt)$. Intuitively, non-decreasable $K$-types are quite close to the trivial one. Figure \ref{u21} illustrates the $U(2, 1)$ case. There are three Weyl chambers $\mathcal{C}^{(i)}$ $(i=0,1,2)$ of $\mathfrak g$ sitting inside the dominant Weyl chamber corresponding to $\Delta^+(\mathfrak k,\mathfrak t)$. Let $(\Delta^+)^{(i)}(\mathfrak g,\mathfrak t)$ be the choice of positive roots corresponding to $\mathcal{C}^{(i)}$. In Figure 1, $\xi_1$ and $\xi_2$ are the two fundamental weights for $(\Delta^+)^{(0)}(\mathfrak{g}, \mathfrak{t})$, and the black dots stand for non-decreasable $K$-types (which turn out to be all u-small), while the white dots stand for other u-small $K$-types. Although there are much fewer non-decreasable $K$-types than u-small ones, the former could be very close to the boundary of the u-small convex hull. Indeed, in Figure 1, the right most black dot can reach the boundary just by adding $\xi_2$.

Now Conjecture 2.1 of \cite{D} reads as: any non-decreasable $K$-type should be u-small. When $G$ is complex or real exceptional, this conjecture has been verified in \cite{D}.
This note aims to prove it for real classical Lie groups. As a result, we have the following.

\begin{thm}\label{thm-main}
Any non-decreasable $K$-type is u-small.
\end{thm}

The paper is organized as follows. Firstly, we handle equal rank groups. Under this assumption, Section \ref{sec-analysis} prepares a series of technical lemmas. Section \ref{sec-equalrank} proves  Theorem \ref{thm-main} for classical equal rank groups. Secondly, Section \ref{sec-nonequalrank} proves Theorem \ref{thm-main} for classical non-equal rank groups.

The lemmas presented in Section \ref{sec-analysis} give delicate analysis of the coefficients of $\gamma$ in terms of the simple roots of $\Delta^+(\mathfrak{g}, \mathfrak{t})$ for real classical Lie algebras, where $\gamma$ is any simple root of $\Delta^+(\mathfrak{k}, \mathfrak{t})$. Though weaker forms meet our need, we deduce their full versions with the hope that they may be helpful elsewhere.

\section{Analysis of compact simple roots}\label{sec-analysis}
We continue with the notation of the introduction.  In this section, we further assume that $G$ is \textbf{equal rank}. That is, $\fra_0=0$.

This section aims to carry out a delicate analysis of the decomposition of the compact simple roots in terms of the simple roots $\Pi$ of $\Delta^+(\mathfrak{g}, \mathfrak{t})$. It will be used later for deducing Theorem \ref{thm-equalrank}.

\subsection{Preliminaries on Vogan diagrams}

\begin{figure}[H]
\centering
\scalebox{0.8}{\includegraphics{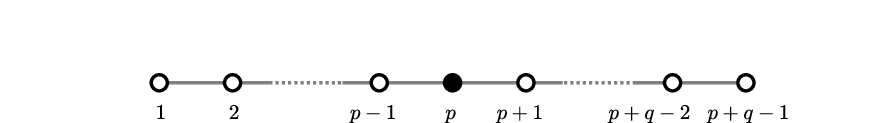}}
\caption{Vogan diagram for $\mathfrak{su}(p, q)$}
\label{upq}
\end{figure}

Since $G$ is assumed to be equal rank, each root $\beta$ in $\Delta(\mathfrak{g}, \mathfrak{t})$ is imaginary. That is, $\beta|_{\fra}$ vanishes. Moreover, $\beta$ is either compact or non-compact. That is, we have either $\frg_{\beta}\subset \frk$ or
$\frg_{\beta}\subset \frp$, respectively. After Knapp \cite{Kn}, the \emph{Vogan diagram} of $\Delta^+(\mathfrak{g}, \mathfrak{t})$ is obtained from the Dynkin diagram of $\Pi$ by painting $\alpha$ to be black whenever it is non-compact. For example, Figure \ref{upq} is the Vogan diagram for $\mathfrak{su}(p, q)$: it paints the $p$-th simple root to be black on the Dynkin diagram of type $A_{n-1}$ with $n=p+q$.

Vogan diagram gives a vivid way to classify non-compact simple real Lie algebras. See Chapter VI of \cite{Kn}. All classical non-compact simple real Lie algebras are listed in Figure 6.1 in \cite{Kn}. We collect several preliminary lemmas.

\begin{lemma}\cite[Lemma 4.5]{HPV}\label{lem.painting.1}
    Let $\Delta^+(\mathfrak k,\mathfrak t)\subset \Delta^+(\mathfrak g,\mathfrak t)$ and let $\Pi$ be the corresponding set of simple roots of $\Delta^+(\mathfrak g,\mathfrak t)$. Let $\alpha\in \Pi$ be non-compact, i.e., a black vertex of the Vogan diagram corresponding to $\Pi$. Then $s_\alpha\Delta^+(\mathfrak g,\mathfrak t)$ is another positive root system containing $\Delta^+(\mathfrak k,\mathfrak t)$, and its set of simple roots is $s_\alpha\Pi$. Furthermore, the Vogan diagram corresponding to $s_\alpha\Pi$ compared with the Vogan diagram corresponding to $\Pi$ has the following:
    \begin{enumerate}
        \item the same colors at the vertex $\alpha$ and at vertices not adjacent to $\alpha$;
        \item opposite colors at vertices adjacent to $\alpha$, unless the neighbor is connected to $\alpha$ by a double edge and longer than $\alpha$;
        \item the same color at the vertex adjacent to $\alpha$ which is connected to $\alpha$ by a double edge and longer than $\alpha$.
    \end{enumerate}
\end{lemma}

\begin{lemma}\cite[Lemma 4.6]{HPV}\label{lem.painting.2}
    Let $(\Delta^+)^{(i)}(\mathfrak g,\mathfrak t)$ and $(\Delta^+)^{(j)}(\mathfrak g,\mathfrak t)$ be positive root systems for $\mathfrak g$ containing $\Delta^+(\mathfrak k,\mathfrak t)$, with corresponding simple roots $\Pi_i$ and $\Pi_j$. Then there are non-compact roots $\alpha_1\in \Pi_i$, $\alpha_2\in s_{\alpha_1}\Pi_i$, $\dots$, $\alpha_k\in s_{\alpha_{k-1}}\cdots s_{\alpha_1}\Pi_i$, such that $(\Delta^+)^{(j)}(\mathfrak g,\mathfrak t)=s_{\alpha_k} s_{\alpha_{k-1}}\cdots s_{\alpha_1}(\Delta^+)^{(i)}(\mathfrak g,\mathfrak t)$.
\end{lemma}

Let $\alpha, \beta\in \Delta^+(\mathfrak g,\mathfrak t)$ such that $\alpha+\beta$ is a root. If exactly one of $\alpha$ and $\beta$ is non-compact, then $\alpha+\beta$ is non-compact; otherwise, $\alpha+\beta$ is compact. This fact is mentioned in (6.99) of \cite{Kn}. The following lemma is a basic application of this fact and we omit the proof.

\begin{lemma}\label{lem.cpt.root.sum}
Let $\gamma$, $\gamma_1$, $\gamma_2$ be roots in $\Delta^+(\mathfrak g, \mathfrak t)$ such that $\gamma=\gamma_1+\gamma_2$. Write $\gamma=\sum_{i}c_i\alpha_i$ as a sum of simple roots with integer coefficients. Then the following statements are equivalent.
    \begin{itemize}
        \item[(a1)] $\gamma$ is compact.
        \item[(a2)] $\gamma_1$ and $\gamma_2$ are both non-compact or both compact.
        \item[(a3)] The sum of all $c_j$ of non-compact $\alpha_j$ is even.
    \end{itemize}
    Correspondingly, the following statements are equivalent.
    \begin{itemize}
        \item[(b1)] $\gamma$ is non-compact.
        \item[(b2)] Exactly one of $\gamma_1$ and $\gamma_2$ is compact.
        \item[(b3)] The sum of all $c_j$ of non-compact $\alpha_j$ is odd.
    \end{itemize}
\end{lemma}

\subsection{Compact simple roots in $\Delta^+(\mathfrak g,\mathfrak t)$}

Appendix C.3 of \cite{Kn} gives an initial choice  $(\Delta^+)^{(0)}(\mathfrak g,\mathfrak t)$ of positive roots of $\Delta(\mathfrak g,\mathfrak t)$  for each simple real Lie algebra $\mathfrak{g}_0$.  At the same time, it offers us an initial choice  $(\Delta^+)^{(0)}(\mathfrak k,\mathfrak t)$ of positive roots of $\Delta(\mathfrak k,\mathfrak t)$. We adopt $(\Delta^+)^{(0)}(\mathfrak k,\mathfrak t)$ as our fixed $\Delta^+(\mathfrak k,\mathfrak t)$, and \textbf{keep the labeling of simple roots} of any $\Delta^+(\mathfrak g,\mathfrak t)$ containing the fixed $\Delta^+(\mathfrak k,\mathfrak t)$ as that of the underlying Dynkin diagram.

\begin{lemma}\label{lem.component.of.gamma}
Let $\gamma$ be a simple root of $\Delta^+(\mathfrak k,\mathfrak t)$. Suppose $\{\alpha_l,\alpha_{l+1},\dots,\alpha_m\}$ are simple roots of $\Delta^+(\mathfrak g,\mathfrak t)$ such that for every integer $k$ $(l\leqslant k\leqslant m)$,
    \begin{itemize}
        \item the sum $\alpha_l+\alpha_{l+1}+\cdots+\alpha_k$ is a root of $\Delta^+(\mathfrak g, \mathfrak t)$;
        \item and $\gamma-(\alpha_l+\alpha_{l+1}+\cdots+\alpha_k)$ is a root of $\Delta^+(\mathfrak g, \mathfrak t)$.
    \end{itemize}
    Then $\alpha_l$ is non-compact, and $\alpha_{l+1}, \dots, \alpha_m$ are all compact. In particular, if $\alpha$ is a simple root of $\Delta^+(\mathfrak g, \mathfrak t)$ such that $\gamma-\alpha$ is a root of $\Delta^+(\mathfrak g, \mathfrak t)$, then $\alpha$ is non-compact.

    We call the set $\{\alpha_l,\alpha_{l+1},\dots,\alpha_m\}$ an \textbf{$\alpha_l$-component of $\gamma$}.
    \begin{proof}
        By Lemma \ref{lem.cpt.root.sum}(a), $(\alpha_l+\alpha_{l+1}+\cdots+\alpha_k)$ and $\gamma-(\alpha_l+\alpha_{l+1}+\cdots+\alpha_k)$ are both non-compact or both compact. Because $\gamma$ is a simple root of $\Delta^+(\mathfrak k,\mathfrak t)$, it cannot be written as a sum of two compact roots. Therefore, for every integer $k$, $\alpha_l+\alpha_{l+1}+\cdots+\alpha_k$ is non-compact. In particular, $\alpha_l$ is non-compact. Notice $\alpha_l+\alpha_{l+1}+\cdots+\alpha_k$ is a sum of two roots $(\alpha_l+\alpha_{l+1}+\cdots+\alpha_{k-1})$ and $\alpha_k$ when $l<k\leqslant m$. By Lemma \ref{lem.cpt.root.sum}(b), $\alpha_k$ is compact for all $l<k\leqslant m$.
    \end{proof}
\end{lemma}

\begin{lemma}\label{lem.cpt.sim.root.upq}
	Let $\mathfrak{g}_0$ be $\mathfrak{su}(p,q)$ with $p+q-1=n$, and $\gamma$ be a simple root of $\Delta^+(\mathfrak k,\mathfrak t)$. Then one of the followings is ture.
	\begin{itemize}
		\item $\gamma$ is a simple root of $\Delta^+(\mathfrak g,\mathfrak t)$.
		\item $\gamma=\alpha_l+\alpha_{l+1}+\cdots+\alpha_m$, where $\alpha_l$ and $\alpha_m$ are non-compact, and the rest summands are compact simple roots.
	\end{itemize}
	\begin{proof}
		Suppose $\gamma$ is not a simple root of $\Delta^+(\mathfrak g,\mathfrak t)$. Then $\gamma=\alpha_l+\alpha_{l+1}+\cdots+\alpha_m$ for two positive integers $l<m$. Notice that $\{\alpha_l,\alpha_{l+1},\dots,\alpha_{m-1} \}$ is an $\alpha_l$-component of $\gamma$, and that $\{\alpha_m,\alpha_{m-1},\dots,\alpha_{l+1}\}$ is an $\alpha_m$-component of $\gamma$. The proof is done by Lemma \ref{lem.component.of.gamma}.
	\end{proof}
\end{lemma}

Next, we are going to describe the simple roots of $\Delta^+(\mathfrak k,\mathfrak t)$ when $\mathfrak{g}_0$ is of type B. As preparation, we list the positive roots $\Delta^+(\mathfrak g,\mathfrak t)$ using the following Dynkin diagram.

\begin{figure}[H]
\centering
\scalebox{0.6}{\includegraphics{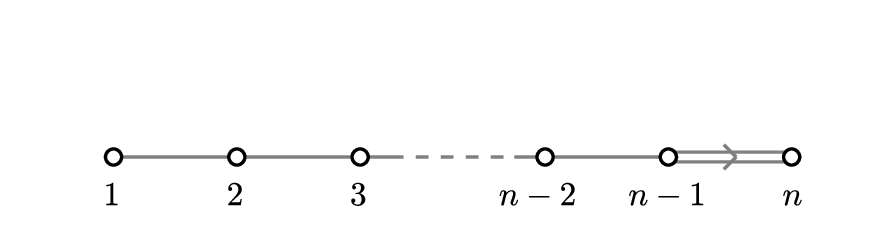}}
\caption{Dynkin diagram for type $B_n$}
\label{Bn}
\end{figure}

The positive roots with respect to this Dynkin diagram consists of the following sets
\begin{itemize}
	\item[(B0)] simple roots $\{\alpha_i|~1\leqslant i\leqslant n\}$;
	\item[(B1)] $\{ \alpha_l+\alpha_{l+1}+\cdots+\alpha_m |~1\leqslant  l\leqslant m< n \}$;
        \item[(B2)] $\{ \alpha_l+\alpha_{l+1}+\cdots+\alpha_n |~1\leqslant  l< n \}$;
	\item[(B3)] $\{\alpha_l+\cdots+\alpha_{m-1}+2(\alpha_m+\cdots+\alpha_n)|~1\leqslant l<m\leqslant n\}$;
\end{itemize}
The short roots are $\alpha_n$ and $\alpha_l+\alpha_{l+1}+\cdots+\alpha_n~(1\leqslant  l< n)$.

\begin{lemma}\label{lem.cpt.sim.root.so2n1}
    Let $\mathfrak{g}_0$ be $\mathfrak{so}(2n,1)$ with $n>1$, and $\gamma$ be a simple root of $\Delta^+(\mathfrak k,\mathfrak t)$. Then one of the followings is true.
    \begin{enumerate}
		\item $\gamma$ is a simple root of $\Delta^+(\mathfrak g,\mathfrak t)$.
		
		\item $\gamma=\alpha_{n-1}+2\alpha_n$, where $\alpha_n$ is non-compact and $\alpha_{n-1}$ is compact.
	\end{enumerate}
    \begin{proof}
    Let $\{\alpha_i^{(0)}\}$ be the simple roots of $(\Delta^+)^{(0)}(\mathfrak g,\mathfrak t)$. Notice that only $\alpha^{(0)}_n$ is non-compact in the set of simple roots. We have listed all the positive roots of a Type B root system. Then Lemma \ref{lem.cpt.root.sum} shows that all short roots are non-compact, and all long roots are compact.

    Now we come back to $\Delta^+(\mathfrak g,\mathfrak t)$. We know three facts: $\gamma$ must be a long root; $\alpha_i$ is compact when $i<n$; and $\alpha_n$ is non-compact. It is clear that $\gamma$ is not in (B2).

    Suppose $\gamma$ is in (B1). Notice that $\{\alpha_l,\alpha_{l+1}, \dots, \alpha_{m-1} \}$ is an $\alpha_l$-component of $\gamma$. By Lemma \ref{lem.component.of.gamma}, $\alpha_l$ is non-compact. This is a contradiction. Thus, $\gamma$ cannot be in (B1).

    Suppose $\gamma$ is in (B3). Because $\{\alpha_m,\alpha_{m+1},\dots,\alpha_n \}$ is an $\alpha_m$-component of $\gamma$, $\alpha_m$ is non-compact. Hence, $m=n$, and now $\gamma=\alpha_l+\cdots+\alpha_{n-1}+2\alpha_n$. When $l<n-1$, $\{\alpha_l,\alpha_{l+1},\dots,\alpha_{n-2} \}$ is an $\alpha_l$-component of $\gamma$. By Lemma \ref{lem.component.of.gamma}, $\alpha_l$ is non-compact. This is a contradiction. When $l=n-1$, we have $\gamma=\alpha_{n-1}+2\alpha_n$. Since $\alpha_n$ is non-compact, $\alpha_{n-1}$ is compact by Lemma \ref{lem.cpt.root.sum}(a).
    \end{proof}
\end{lemma}

\begin{lemma}\label{lem.cpt.sim.root.so2p2q+1}
    Let $\mathfrak{g}_0$ be $\mathfrak{so}(2p,2q+1)$ with $p+q=n>1$, and $\gamma$ be a simple root of $\Delta^+(\mathfrak k,\mathfrak t)$. Then one of the followings is true.
    \begin{enumerate}
		\item $\gamma$ is a simple root of $\Delta^+(\mathfrak g,\mathfrak t)$.
		\item $\gamma=\alpha_l+\alpha_{l+1}+\cdots+\alpha_m$, where $\alpha_l$ and $\alpha_m$ are non-compact, and the rest summands are compact.
        \item $\gamma=\alpha_l+\cdots+\alpha_{m-1}+2(\alpha_m+\cdots+\alpha_n)$, where $\alpha_l$, $\alpha_{m-1}$ and $\alpha_m$ are non-compact, and the rest summands are compact.
        \item $\gamma=\alpha_{m-1}+2(\alpha_m+\cdots+\alpha_n)$, where $\alpha_m$ is non-compact, and the rest summands are compact.
	\end{enumerate}
    \begin{proof}
    Suppose $\gamma$ is in (B1) or (B2). The proof is same to the case in Lemma \ref{lem.cpt.sim.root.upq}.

    Suppose $\gamma$ is in (B3). Notice that $\{\alpha_m,\alpha_{m+1},\dots,\alpha_n \}$ is an $\alpha_m$-component of $\gamma$. By Lemma \ref{lem.component.of.gamma}, $\alpha_m$ is non-compact and all $\alpha_{m+1},\dots,\alpha_n$ are compact. When $l<m-1$, $\{\alpha_l,\alpha_{l+1},\dots,\alpha_{m-2} \}$ is an $\alpha_l$-component of $\gamma$. By Lemma \ref{lem.component.of.gamma}, $\alpha_l$ is non-compact and all $\alpha_{l+1},\dots,\alpha_{m-2}$ are compact. By Lemma \ref{lem.cpt.root.sum}(a), $\alpha_{m-1}$ is non-compact. When $l=m-1$, we have $\gamma=\alpha_{m-1}+2(\alpha_m+\cdots+\alpha_n)$. Note that $\alpha_{m-1}$ is compact by (a) Lemma \ref{lem.cpt.root.sum}.
    \end{proof}
\end{lemma}

Next, we are going to describe the simple roots of $\Delta^+(\mathfrak k,\mathfrak t)$ when $\mathfrak{g}_0$ is of type C. As preparation, we list the positive roots $\Delta^+(\mathfrak g,\mathfrak t)$ using the following Dynkin diagram.

\begin{figure}[H]
\centering
\scalebox{0.6}{\includegraphics{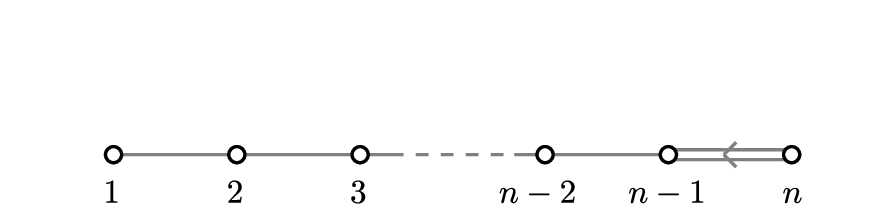}}
\caption{Dynkin diagram for type $C_n$}
\label{Cn}
\end{figure}

The positive roots with respect to this Dynkin diagram are
\begin{itemize}
	\item[(C0)] simple roots $\{\alpha_i|1\leqslant i\leqslant n\}$;
	\item[(C1)] $\{ \alpha_l+\alpha_{l+1}+\cdots+\alpha_m |1\leqslant  l< m\leqslant n \}$;
	\item[(C2)] $\{\alpha_l+\cdots+\alpha_{m-1}+2(\alpha_m+\cdots+\alpha_{n-1})+\alpha_n|1\leqslant  l< m< n \}$;
	\item[(C3)] $\{2(\alpha_m+\cdots+\alpha_{n-1})+\alpha_n|1\leqslant m< n\}$.
\end{itemize}
The long roots are $\alpha_n$ and $2(\alpha_m+\cdots+\alpha_{n-1})+\alpha_n~(1\leqslant m< n)$.

\begin{lemma}\label{lem.cpt.sim.root.sp2nr}
	Let $\mathfrak{g}_0$ be $\mathfrak{sp}(n,\mathbb R)$ with $n>2$, and $\gamma$ be a simple root of $\Delta^+(\mathfrak k,\mathfrak t)$. Then one of the followings is true.
	\begin{enumerate}
		\item $\gamma$ is a simple root of $\Delta^+(\mathfrak g,\mathfrak t)$.
		
		\item $\gamma=\alpha_l+\alpha_{l+1}+\cdots+\alpha_m$, where $\alpha_l$ and $\alpha_m$ are non-compact, and the rest summands are compact.
	\end{enumerate}
	\begin{proof}
        Let $\{\alpha_i^{(0)}\}$ be the simple roots of $(\Delta^+)^{(0)}(\mathfrak g,\mathfrak t)$. Notice that only $\alpha^{(0)}_n$ is non-compact in the set of simple roots. All the long positive roots are $\alpha_n^{(0)}$ and $2(\alpha_m^{(0)}+\cdots+\alpha_{n-1}^{(0)})+\alpha_n^{(0)}~(1\leqslant m< n)$. By Lemma \ref{lem.cpt.root.sum}(b), they are all non-compact.

        Now we come back to $\Delta^+(\mathfrak g,\mathfrak t)$. We know two facts: $\gamma$ must be a short root; and $\alpha_n$ is non-compact. It is clear that $\gamma$ is not in (C3).

        Suppose $\gamma$ is in (C1), then the proof is same to the case in Lemma \ref{lem.cpt.sim.root.upq}.

        Suppose $\gamma$ is in (C2). Because $\{\alpha_l,\alpha_{l+1},\dots,\alpha_{n-1}\}$ is an $\alpha_l$-component of $\gamma$, $\alpha_m$ is compact. Notice that $\gamma=(\gamma-\alpha_m)+\alpha_m$ is a sum of two compact roots. This contradicts to the fact that $\gamma$ is a simple root of $\Delta^+(\mathfrak k,\mathfrak t)$. Hence, $\gamma$ cannot be in (C2).
	\end{proof}
\end{lemma}

\begin{lemma}\label{lem.cpt.sim.root.sppq}
		Let $\mathfrak{g}_0$ be $\mathfrak{sp}(p,q)$ with $p+q=n>2$, and $\gamma$ be a simple root of $\Delta^+(\mathfrak k,\mathfrak t)$. Then one of the followings is true.
	\begin{enumerate}
		\item $\gamma$ is a simple root of $\Delta^+(\mathfrak g,\mathfrak t)$.
		\item $\gamma=\alpha_l+\alpha_{l+1}+\cdots+\alpha_m$, where $\alpha_l$ and $\alpha_m$ are non-compact, and the rest summands are compact.
		\item  $\gamma=2(\alpha_l+\alpha_{l+1}+\cdots+\alpha_{n-1})+\alpha_n$, where $\alpha_l$ is non-compact, and the rest summands are compact.
	\end{enumerate}
	\begin{proof}
        Suppose $\gamma$ is in (C1), then the proof is same to the case in Lemma \ref{lem.cpt.sim.root.upq}.

		Same to the case of $\mathfrak{sp}(n,\mathbb R)$, $\gamma$ cannot be in (C2).

        Suppose $\gamma$ is in (C3). Notice that $\{\alpha_m,\alpha_{m+1},\dots,\alpha_n\}$ is an $\alpha_m$-component of $\gamma$. By Lemma \ref{lem.component.of.gamma}, $\alpha_m$ is non-compact and all $\alpha_{m+1},\dots,\alpha_n$ are compact.
	\end{proof}
\end{lemma}

Next, we are going to describe the simple roots of $\Delta^+(\mathfrak k,\mathfrak t)$ when $\mathfrak{g}_0$ is of type D. As preparation, we list the positive roots $\Delta^+(\mathfrak g,\mathfrak t)$ using the following Dynkin diagram.

\begin{figure}[H]
\centering
\scalebox{0.6}{\includegraphics{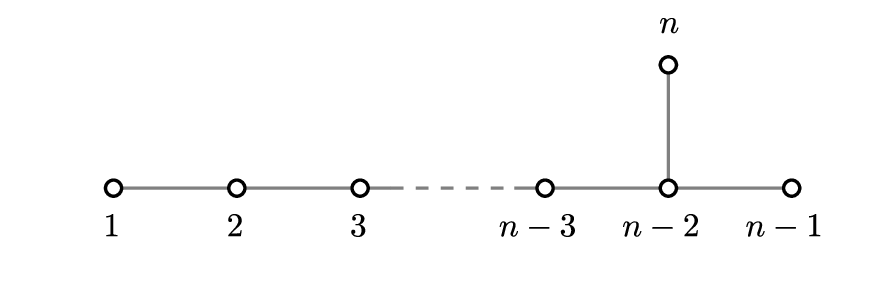}}
\caption{Dykin diagram for type $D_n$}
\label{Dn}
\end{figure}

The positive roots with respect to this Dynkin diagram are
\begin{itemize}
    \item[(D0)] simple roots $\{\alpha_i|~1\leqslant i\leqslant n\}$;
    \item[(D1)] $\{\alpha_l+\alpha_{l+1}+\cdots+\alpha_m|~1\leqslant l<m\leqslant n-1\}$;
    \item[(D2)] $\{\alpha_l+\alpha_{l+1}+\cdots+\alpha_{n-2}+\alpha_n |~l<n-1 \}$;
    \item[(D3)] $\{ \alpha_l+\alpha_{l+1}+\cdots+\alpha_n|~l<n-1 \}$;
    \item[(D4)] $\{\alpha_l+\cdots+\alpha_{m-1}+2(\alpha_m+\cdots+\alpha_{n-2})+\alpha_{n-1}+\alpha_n |~l<m<n-1\}$.
\end{itemize}

\begin{lemma}\label{lem.cpt.sim.root.so2n}
    Let $\mathfrak{g}_0$ be $\mathfrak{so}^\ast(2n)$ with $n>3$, and $\gamma$ be a simple root of $\Delta^+(\mathfrak k,\mathfrak t)$. Then one of the followings is true.
    \begin{enumerate}
        \item $\gamma$ is a simple root of $\Delta^+(\mathfrak g,\mathfrak t)$.
		
		\item $\gamma=\alpha_l+\alpha_{l+1}+\cdots+\alpha_m~(l<m<n)$, where $\alpha_l$ and $\alpha_m$ are non-compact, and the rest summands are compact.
        \item $\gamma=\alpha_l+\alpha_{l+1}+\cdots+\alpha_{n-2}+\alpha_n~(l<n-1)$, where $\alpha_l$ and $\alpha_n$ are non-compact, and the rest summands are compact.
    \end{enumerate}
    \begin{proof}
        Let $\Pi_0=\{\alpha_i^{(0)}\}$ be the simple roots of $(\Delta^+)^{(0)}(\mathfrak g,\mathfrak t)$. Notice that only $\alpha^{(0)}_n$ is non-compact in the set of simple roots. Let $\Pi$ be the set of simple roots of $\Delta^+(\mathfrak g,\mathfrak t)$. Lemmas \ref{lem.painting.1} and \ref{lem.painting.2} give an algorithm to obtain the Vogan diagram for $\Pi$ from that of $\Pi_0$. Notice that only the reflection corresponding to the simple root at the $(n-2)$-th vertex can change the color of the $(n-1)$-th and the $n$-th vertices. When we proceed with such a reflection, the color of the $(n-1)$-th and the $n$-th vertices will change simultaneously. Because in the Vogan diagram of the initial $\Pi_0$ has only one painted vertex which is the one for $\alpha_n^{(0)}$, the vertices for $\alpha_n$ and $\alpha_{n-1}$ must have opposite colors. Therefore, exactly one of $\alpha_{n-1}$ and $\alpha_n$ is non-compact.

        Suppose $\gamma$ is in (D1) or (D2), then the proof is same to the case in Lemma \ref{lem.cpt.sim.root.upq}.

        Suppose $\gamma$ is in (D3). Notice that both $\gamma-\alpha_{n-1}$ and $\gamma-\alpha_n$ are positive roots. Then both $\alpha_{n-1}$ and $\alpha_n$ have to be non-compact by Lemma \ref{lem.component.of.gamma}. Hence, we have a contradiction and $\gamma$ cannot be in (D3).

        Suppose $\gamma$ is in (D4). We have proved that exactly one of $\alpha_{n-1}$ and $\alpha_n$ is non-compact. Without loss of generality, we may assume $\alpha_{n-1}$ is non-compact. Since $\{\alpha_m,\alpha_{m+1},\dots,\alpha_{n-1}\}$ is an $\alpha_m$-component of $\gamma$, the root $\alpha_{n-1}$ should be compact. Hence, we have a contradiction and $\gamma$ cannot be in (D4).
    \end{proof}
\end{lemma}

\begin{lemma}\label{lem.cpt.sim.root.so2p2q}
    Let $\mathfrak{g}_0$ be $\mathfrak{so}(2p,2q)$ with $p+q=n>3$, and $\gamma$ be a simple root of $\Delta^+(\mathfrak k,\mathfrak t)$. Then one of the followings is true.
    \begin{enumerate}
        \item $\gamma$ is a simple root of $\Delta^+(\mathfrak g,\mathfrak t)$.
		
		\item $\gamma=\alpha_l+\alpha_{l+1}+\cdots+\alpha_m~(l<m<n)$, where $\alpha_l$ and $\alpha_m$ are non-compact, and the rest summands are compact.

        \item $\gamma=\alpha_l+\alpha_{l+1}+\cdots+\alpha_{n-2}+\alpha_n~(l<n-1)$, where $\alpha_l$ and $\alpha_n$ are non-compact, and the rest summands are compact.

        \item $\gamma=\alpha_l+\alpha_{l+1}+\cdots+\alpha_n~(l<n-1)$, where $\alpha_l$, $\alpha_{n-2}$, $\alpha_{n-1}$ and $\alpha_n$ are non-compact, and the rest summands are compact.

        \item $\gamma=\alpha_{n-2}+\alpha_{n-1}+\alpha_n$, where $\alpha_{n-1}$ and $\alpha_n$ are non-compact, and $\alpha_{n-2}$ is compact.

        \item $\gamma=\alpha_l+\alpha_{l+1}+\cdots+\alpha_{m-1}+2(\alpha_{m}+\cdots+\alpha_{n-2})+\alpha_{n-1}+\alpha_n~(l<m<n-1)$, where $\alpha_{l}$, $\alpha_{m-1}$ and $\alpha_m$ are non-compact, and the rest summands are compact.

        \item $\gamma=\alpha_l+2(\alpha_{l+1}+\cdots+\alpha_{n-2})+\alpha_{n-1}+\alpha_n~(l<n-2)$, where only $\alpha_{l+1}$ is non-compact, and the rest summands are compact.
    \end{enumerate}
    \begin{proof}
        Using the same method in Lemma \ref{lem.cpt.sim.root.so2n}, we can prove that $\alpha_{n-1}$ and $\alpha_n$ are both non-compact or both compact.

        Suppose $\gamma$ is in (D1) or (D2), then the proof is same to the case in Lemma \ref{lem.cpt.sim.root.upq}.

        Suppose $\gamma$ is in (D3). Notice both $\gamma-\alpha_{n-1}$ and $\gamma-\alpha_n$ are positive roots. Then both $\alpha_{n-1}$ and $\alpha_n$ has to be non-compact by Lemma \ref{lem.component.of.gamma}. When $l<n-2$, $\{\alpha_l,\alpha_{l+1},\dots,\alpha_{n-3}\}$ is an $\alpha_l$-component of $\alpha_l$-component of $\gamma$. So, $\alpha_l$ is non-compact, and all $\alpha_{l+1},\dots,\alpha_{n-3}$ are compact. By Lemma \ref{lem.cpt.root.sum}(a), $\alpha_{n-2}$ is non-compact. This gives item (4). When $l=n-2$, $\gamma=\alpha_{n-2}+\alpha_{n-1}+\alpha_n$, by Lemma \ref{lem.cpt.root.sum}(a), $\alpha_{n-2}$ is compact. This gives item (5).

        Suppose $\gamma$ is in (D4). Because $\{\alpha_m,\alpha_{m+1},\dots,\alpha_{n-2},\alpha_{n-1}\}$ and $\{\alpha_m,\alpha_{m+1},\dots,\alpha_{n-2},\alpha_n\}$ are two $\alpha_m$-components of $\gamma$, we know $\alpha_m$ is non-compact, and all $\alpha_{m+1},\dots,\alpha_{n-1},\alpha_n$ are compact. When $l<m-1$, $\{\alpha_l,\alpha_{l+1},\dots,\alpha_{m-2} \}$ is an $\alpha_l$-component of $\gamma$. So, $\alpha_l$ is non-compact, and  $\alpha_{l+1},\dots,\alpha_{m-2}$ are all compact. By Lemma \ref{lem.cpt.root.sum}(a), $\alpha_{m-1}$ is non-compact. This gives item (6). When $l=m-1$, By Lemma \ref{lem.cpt.root.sum}(a), $\alpha_l$ is compact. This gives item (7).
    \end{proof}
\end{lemma}

\section{Equal rank groups}
\label{sec-equalrank}

 We continue to assume that $G$ is equal rank. Write the set of simple roots of $\Delta^+(\mathfrak g,\mathfrak t)$ as $\Pi=\{\alpha_1,\dots,\alpha_{n}\}$, where $n$ is the rank of $G$. Let $\{\xi_1,\xi_2,\dots,\xi_n\}$ be the fundamental weights in the sense that
\begin{equation*}
	\frac{2(\xi_i,\alpha_j)}{(\alpha_j,\alpha_j)}=\begin{cases}
		1& i=j,\\
		0&i\neq j.
	\end{cases}
\end{equation*}
We say a weight $\mu$ is $\mathfrak g$-dominant (resp. $\mathfrak k$-dominant) if $\mu$ is dominant with respect to $\Delta^+(\mathfrak g,\mathfrak t)$ (resp. $\Delta^+(\mathfrak k,\mathfrak t)$).
Now we are able to distill the following lemma based on the results of the previous section.

\begin{lemma}\label{lem.cpt.root.coef}
	Let $\mathfrak{g}_0$ be the Lie algebra of $G$. Take a simple root $\gamma$ of $\Delta^+(\mathfrak k,\mathfrak t)$ and write it as a sum $\gamma=\sum_{j}c_j\alpha_j$. Then
	\begin{itemize}
		\item[(a)]
		$c_i(\alpha_i,\alpha_i)\geqslant (\gamma,\gamma)$ when $c_i\neq 0$.
		\item[(b)]
		$\frac{2(\xi_i,\gamma)}{(\gamma,\gamma)}$ is $2$, $1$ or $0$.
		\item[(c)] Let $\mu$ be an integral weight of $\mathfrak g$, then $\frac{2(\mu,\gamma)}{(\gamma,\gamma)}\in\mathbb Z$.
	\end{itemize}

	\begin{proof}
		For (a), we may assume $\gamma$ is not a simple root of $\Delta^+(\mathfrak g,\mathfrak t)$. Since $c_i$ is a positive integer, it suffices to prove the case when $\gamma$ is a long root and $\alpha_i$ is a short root. This case happens in Lemma \ref{lem.cpt.sim.root.so2n1}(2), Lemma \ref{lem.cpt.sim.root.so2p2q+1}(3)(4) and Lemma \ref{lem.cpt.sim.root.sppq}(3). For $\mathfrak{so}(2n,1)$ and $\mathfrak{so}(2p,2q+1)$, where $\alpha_n$ is short and we have $c_n=2$ and $2(\alpha_n,\alpha_n)=(\gamma,\gamma)$. For $\mathfrak{sp}(p,q)$, where $\alpha_n$ is long and we have $c_i=2$ and $2(\alpha_i,\alpha_i)=(\gamma,\gamma)$ for $i\neq n$. Then (a) follows.
		
		For (b), it is sufficient to discuss the case when $c_i\neq 0$. We compute
		$$
		\frac{2(\xi_i,\gamma)}{(\gamma,\gamma)}=\frac{2(\xi_i,\alpha_i)}{(\alpha_i,\alpha_i)}\cdot c_i\cdot \frac{(\alpha_i,\alpha_i)}{(\gamma,\gamma)}= c_i\cdot \frac{(\alpha_i,\alpha_i)}{(\gamma,\gamma)}.
		$$
		There are only four situations.
		\begin{itemize}
		\item $(\alpha_i,\alpha_i)=(\gamma,\gamma)$ and $c_i=1$;
		\item only in $\mathfrak{so}(2n,1)$, $\mathfrak{so}(2p,2q+1)$ and $\mathfrak{sp}(p,q)$, we could have $2(\alpha_i,\alpha_i)=(\gamma,\gamma)$ and $c_i=2$;
        \item only in $\mathfrak{so}(2p,2q+1)$ and $\mathfrak{so}(2p,2q)$, we could have $(\alpha_i,\alpha_i)=(\gamma,\gamma)$ and $c_i=2$.
        \item only in $\mathfrak{so}(2p,2q+1)$, $\mathfrak{sp}(n,\mathbb R)$ and $\mathfrak{sp}(p,q)$, we could have $(\alpha_i,\alpha_i)=2(\gamma,\gamma)$ and $c_i=1$.
		\end{itemize}
        In the first two cases, $\frac{2(\xi_i,\gamma)}{(\gamma,\gamma)}=1$; and in the last two cases, $\frac{2(\xi_i,\gamma)}{(\gamma,\gamma)}=2$.
		
		(c) is a direct corollary of (b) since $\mu$ is a sum of fundamental weights with integer coefficients.
	\end{proof}
\end{lemma}

Now let us deduce the following result.

\begin{thm}\label{thm-equalrank}
Assume that $G$ is classical and equal rank. Then any non-decreasable $K$-type is u-small.
	\begin{proof}
		Let $\mu$ be a non-decreasable integral weight. By definition, $\mu$ is $\mathfrak k$-dominant, and $(\mu,\alpha)$ is an integer for any root $\alpha$. Let $\Delta^+(\mathfrak g,\mathfrak t)$ be a choice of positive roots such that $\mu+2\rho_{\rm c}$ is $\mathfrak g$-dominant. Let $\{\alpha_i\}$ be the simple roots of $\Delta^+(\mathfrak g,\mathfrak t)$ and $\{\xi_j\}$ be the corresponding fundamental roots. Then for any fundamental weight $\xi_i$, either $\mu-\xi_i$ is not $\mathfrak k$-dominant or $\mu-\xi_i+2\rho_{\rm c}$ is not $\mathfrak g$-dominant. Let us proceed according to the following cases.
		
		\noindent Case I.  The simple root $\alpha_i$ is compact. Then $(2\rho_{\rm c},\alpha_i)=(\alpha_i,\alpha_i)$. We compute
		\begin{equation*}
			(\mu-\xi_i+2\rho_{\rm c},\alpha_j)=
			\begin{cases}
				(\mu+2\rho_{\rm c},\alpha_j)\geqslant 0, & j\neq i,\\
				(\mu,\alpha_i)+\frac{(\alpha_i,\alpha_i)}{2}> 0, & j=i.
			\end{cases}
		\end{equation*}
		Hence, $\mu-\xi_i+2\rho_{\rm c}$ must be $\mathfrak g$-dominant. Therefore, $\mu-\xi_i$ is not $\mathfrak k$-dominant.
        As a consequence, $(\mu-\xi_i,\gamma)<0$ for some simple root $\gamma$ of $\Delta^+(\mathfrak k,\mathfrak t)$.

        Notice that the proof henceforth actually does not require $\alpha_i$ to be compact.
        By items (b) and (c) of Lemma \ref{lem.cpt.root.coef}, $\frac{2(\xi_i,\gamma)}{(\gamma,\gamma)}$ is $2$, $1$ or $0$; moreover, $\frac{2(\mu,\gamma)}{(\gamma,\gamma)}$ is a non-negative integer. Hence, we have two sub-cases.
        \begin{itemize}
            \item[(a)] $(\mu,\gamma)=0$ and $(\xi_i,\gamma)>0$. Write $\gamma=\sum_ic_i\alpha_i$, where $c_i$ are non-negative integers. Then
$$
		0=(\mu,\gamma)=(\mu+2\rho_{\rm c},\gamma)-(2\rho_{\rm c},\gamma)=\sum_j(\mu+2\rho_{\rm c},c_j\alpha_j)-(2\rho_{\rm c},\gamma).
$$
		Since $(2\rho_{\rm c},\gamma)=(\gamma,\gamma)$ and $(\mu+2\rho_{\rm c},\alpha_j)\geqslant 0$ for all $j$, we have
        $$
        c_i(\mu+2\rho_{\rm c},\alpha_i)\leqslant (\gamma,\gamma).
        $$
        Notice that $(2\rho,\alpha_i)=(\alpha_i,\alpha_i)$ and $\rho=\rho_{\rm n}+\rho_{\rm c}$. Therefore,
		\begin{equation*}
			(2\rho_{\rm n}-\mu,\alpha_i)=\big(2\rho-(\mu+2\rho_{\rm c}),\alpha_i\big)=(\alpha_i,\alpha_i)-(\mu+2\rho_{\rm c},\alpha_i)\geqslant (\alpha_i,\alpha_i)-c_i^{-1}(\gamma,\gamma).
		\end{equation*}
		By Lemma \ref{lem.cpt.root.coef}(a), $(2\rho_{\rm n}-\mu,\alpha_i)\geqslant 0$.

        \item[(b)] $\frac{2(\mu,\gamma)}{(\gamma,\gamma)}=1$ and $\frac{2(\xi_i,\gamma)}{(\gamma,\gamma)}=2$. By Lemma \ref{lem.cpt.root.coef}(b), we have $\frac{(\gamma,\gamma)}{c_i}=\frac{(\alpha_i,\alpha_i)}{2}$. Then
$$
\frac{(\gamma,\gamma)}{2}=(\mu,\gamma)=(\mu+2\rho_{\rm c},\gamma)-(2\rho_{\rm c},\gamma)=\sum_j(\mu+2\rho_{\rm c},c_j\alpha_j)-(2\rho_{\rm c},\gamma).
$$
		Since $(2\rho_{\rm c},\gamma)=(\gamma,\gamma)$ and $(\mu+2\rho_{\rm c},\alpha_j)\geqslant 0$ for all $j$, we have
        $$
        c_i(\mu+2\rho_{\rm c},\alpha_i)\leqslant \frac{3}{2}(\gamma,\gamma).
        $$
        Because $\frac{(\gamma,\gamma)}{c_i}=\frac{(\alpha_i,\alpha_i)}{2}$, we have
        $$
        (\mu+2\rho_{\rm c},\alpha_i)\leqslant \frac{3}{4}(\alpha_i,\alpha_i).
        $$ Therefore,
		\begin{equation*}
			(2\rho_{\rm n}-\mu,\alpha_i)=(\alpha_i,\alpha_i)-(\mu+2\rho_{\rm c},\alpha_i)\geqslant \frac{(\alpha_i,\alpha_i)}{4}>0.
		\end{equation*}
        \end{itemize}

		\noindent Case II. The simple root $\alpha_i$ is non-compact. We move on according to the following sub-cases:
		\begin{itemize}
			\item [(a)] $\mu-\xi_i+2\rho_{\rm c}$ is not $\mathfrak g$-dominant. We must have $(\mu+2\rho_{\rm c},\alpha_i)=0$. Then
			\begin{equation*}
				(2\rho_{\rm n}-\mu,\alpha_i)=\big(2\rho-(\mu+2\rho_{\rm c}),\alpha_i\big)=(\alpha_i,\alpha_i)\geqslant 0.
			\end{equation*}
			
			\item [(b)] $\mu-\xi_i$ is not $\mathfrak k$-dominant, i.e.,  there exists a simple root $\gamma$ of $\Delta^+(\mathfrak k,\mathfrak t)$ such that  $(\mu-\xi_i,\gamma)<0$. Then similar to Case I, one deduces that $(2\rho_{\rm n}-\mu,\alpha_i)\geqslant 0$.
		\end{itemize}
		To sum up, we always have $(2\rho_{\rm n}-\mu,\alpha_i)\geqslant 0$ for all $i$, meaning that $2\rho_{\rm n}-\mu$ is $\mathfrak g$-dominant. It follows from Theorem 6.7(e) of \cite{SV} that $\mu$ is unitarily small.
	\end{proof}
\end{thm}

\section{Non-equal rank groups}
\label{sec-nonequalrank}
Assuming that $G$ is \emph{not} equal rank, this section aims to prove the following.

\begin{thm}\label{thm-nonequalrank}
Assume that $G$ is classical and not equal rank. Then any non-decreasable $K$-type is u-small.
\end{thm}

The above theorem will be proven case by case.

\subsection{$\mathfrak{g}_0=\mathfrak{sl}(2n+1, \bbR)$}
As on page 460 of \cite{Vog86}, we  arrange that
$$
\Delta^+(\mathfrak{k}, \mathfrak{t})=\{e_i \pm e_j \mid 1\leqslant i<j\leqslant n\}\cup\{e_i\mid 1 \leqslant i\leqslant n\},$$
and that
$$\Delta^+(\mathfrak{p}, \mathfrak{t})=\Delta^+(\mathfrak{k}, \mathfrak{t})\cup\{2e_i\mid 1\leqslant i\leqslant n\}.
$$
Thus $W(\mathfrak{g}, \mathfrak{t})^1=\{e\}$ and $\xi_i=e_1+\cdots+e_i$ for $1\leqslant i\leqslant n-1$, while $\xi_n=\frac{1}{2}(e_1+\cdots+e_n)$. Moreover, we have $$
\rho_{\rm c}=\xi_1+\xi_2+\cdots+\xi_n.
$$
Note that $\xi_1, \dots, \xi_n$ are also the fundamental weights of $\Delta^+(\mathfrak{k}, \mathfrak{t})$.
Thus for any $\mu$ which is $\mathfrak{k}$-dominant integral,  the weight  $\mu+2\rho_{\rm c} -\xi_i$ is $\mathfrak{g}$-dominant for any $1\leqslant i\leqslant n$. Therefore, $\mu$ is non-decreasable if and only if $\mu-\xi_i$ is not $\mathfrak{k}$-dominant for any $1\leqslant i\leqslant n$. Then $\mu$ has to be the zero weight, which is u-small.

\subsection{$\mathfrak{g}_0=\mathfrak{sl}(2n, \bbR)$}
As on page 460 of \cite{Vog86}, we arrange that
$$
\Delta^+(\mathfrak{k}, \mathfrak{t})=\{e_i \pm e_j \mid 1\leqslant i<j\leqslant n\}.$$
The corresponding fundamental weights are $\varpi_i = e_1+\cdots+e_i$ for $1\leqslant i\leqslant n-2$,
$$
\varpi_{n-1}=\frac{1}{2}(e_1+\cdots+e_{n-1}-e_n), \quad \varpi_n=\frac{1}{2}(e_1+\cdots+e_{n-1}+e_n).
$$
We take $\mu=[m_1, m_2, \dots, m_n]$, standing for the weight $m_1\varpi_1+\cdots+m_n\varpi_n$. Here each $m_i$ is a non-negative integer.

Case I. $\Delta^+(\mathfrak{p}, \mathfrak{t})=\Delta^+(\mathfrak{k}, \mathfrak{t})\cup \{2 e_1, \dots, 2 e_{n-1}, 2 e_n\}$. Then $\xi_i=\varpi_i$ for $1\leqslant i\leqslant n-2$, $\xi_{n-1}=\varpi_{n-1}+\varpi_n$ and $\xi_n = 2\varpi_n$. Now
$$
\mu+2\rho_{\rm c}=(m_1+2)\xi_1+\cdots+(m_{n-1}+2)\xi_{n-1}+\frac{m_n-m_{n-1}}{2}\xi_n.
$$
Thus $\mu+2\rho_{\rm c}$ is $\mathfrak{g}$-dominant if and only if $m_{n-1}\leqslant m_n$.
Assume that $\mu$ is non-decreasable. It is evident that $\mu+2\rho_{\rm c}-\xi_i$ is $\mathfrak{g}$-dominant for $1\leqslant i\leqslant n-1$. Thus $\mu-\xi_i$ is not $\mathfrak{k}$-dominant for $1\leqslant i\leqslant n-1$. Therefore, we have that $m_i=0$ for $1\leqslant i\leqslant n-2$ and $\min\{m_{n-1}, m_n\}=0$. To move on, we have the following sub-cases:
\begin{itemize}
\item[(a)] $\mu-\xi_n$ is not $\mathfrak{k}$-dominant. Then $m_n=0$ or $1$. Since $m_{n-1}\leqslant m_n$ and $\min\{m_{n-1}, m_n\}=0$, we only have $(m_{n-1}, m_n)=(0, 0)$ or $(0, 1)$. Thus $\mu$ is trivial or $\varpi_n$. Both of them are u-small.

\item[(b)] $\mu + 2\rho_{\rm c}-\xi_n$ is not $\mathfrak{g}$-dominant. Then $m_n-m_{n-1}<2$. Thus it equals $0$ or $1$. Again, we have $(m_{n-1}, m_n)=(0, 0)$ or $(0, 1)$. That is, $\mu$ is trivial or $\varpi_n$. Both of them are u-small.
\end{itemize}

Case II. $\Delta^+(\mathfrak{p}, \mathfrak{t})=\Delta^+(\mathfrak{k}, \mathfrak{t})\cup \{2 e_1, \dots, 2 e_{n-1}, -2 e_n\}$. Then $\xi_i=\varpi_i$ for $1\leqslant i\leqslant n-2$, $\xi_{n-1}=\varpi_{n-1}+\varpi_n$ and $\xi_n = 2\varpi_{n-1}$. Now
$$
\mu+2\rho_{\rm c}=(m_1+2)\xi_1+\cdots+(m_n+2)\xi_{n-1}+\frac{m_{n-1}-m_n}{2}\xi_n.
$$
Thus $\mu+2\rho_{\rm c}$ is $\mathfrak{g}$-dominant if and only if $m_{n-1}\geq m_n$.
Similar to the previous case, one computes that there are only two non-decreasable $\mathfrak{k}$-types: the trivial one and $\varpi_{n-1}$. Both of them are u-small.

\subsection{$\mathfrak{g}_0=\mathfrak{sl}(n, \bbH)$}
As on page 462 of \cite{Vog86}, we arrange that
$$
\Delta^+(\mathfrak{p}, \mathfrak{t})=\{e_i \pm e_j \mid 1\leqslant i<j\leqslant n\}, \quad \Delta^+(\mathfrak{k}, \mathfrak{t})=\Delta^+(\mathfrak{p}, \mathfrak{t})\cup\{2e_i\mid 1\leqslant i\leqslant n\}.
$$
Thus $W(\mathfrak{g}, \mathfrak{t})^1=\{e\}$ and $\xi_i=e_1+\cdots+e_i$ for $1\leqslant i\leqslant n$. Moreover, we have $$
\rho_{\rm c}=\xi_1+\xi_2+\cdots+\xi_n.
$$
Note that $\xi_1, \dots, \xi_n$ are also the fundamental weights of $\Delta^+(\mathfrak{k}, \mathfrak{t})$.
Similar to the $\mathfrak{sl}(2n+1, \bbR)$ case, only the trivial $\mathfrak{k}$-type is non-decreasable. It is u-small.

\subsection{$\mathfrak{g}_0=\mathfrak{so}(2p+1, 2q+1)$}

By restricting the initial choice of positive roots as in Appendix C.3 of \cite{Kn} to $\frt$, we have
$$
(\Delta^+)^{(0)}(\mathfrak g, \mathfrak t)=\{ e_i\pm e_j|1\leqslant i<j\leqslant p+q\}\cup \{e_k|1\leqslant k\leqslant p+q\}.
$$
The simple roots for $\Delta^+(\mathfrak k, \mathfrak t)$ are $e_p$ (if $p>0$), $e_{p+q}$ and all $e_i-e_{i+1}$ with $1\leqslant i\leqslant p-1$ or $p+1\leqslant i\leqslant p+q-1$.

Notice that $(\Delta^+)^{(0)}(\mathfrak g, \mathfrak t)$ is of type $B_{p+q}$, and $\Delta^+(\mathfrak k,\mathfrak t)\subset (\Delta^+)^{(0)}(\mathfrak g, \mathfrak t)$. The analysis of the simple roots for $\Delta^+(\mathfrak k, \mathfrak t)$ is covered by Lemma \ref{lem.cpt.sim.root.so2p2q+1}. Then, similar to Theorem \ref{thm-equalrank},  one proves Theorem \ref{thm-nonequalrank} for $\mathfrak{so}(2p+1, 2q+1)$.

\section*{Funding}
Dong is supported by the National Natural Science Foundation of China (grant no.~12171344, 2022-2025).



\begin{thebibliography}{99}

\bibitem{D}
C.-P.~Dong, \emph{On the Helgason-Johnson bound}, Israel J. Math \textbf{254} (2023), 373--397.

\bibitem{DDW}
Y.-H.~Ding, C.-P.~Dong and L.~Wei, \emph{Dirac series of $E_{7(7)}$}, Israel J. Math \textbf{265} (2025), 347--378.

\bibitem{DDDLY}
Y.-H. Ding, C.-P.~Dong, C.~Du, Y.-Z.~Luan and L.~Yang, \emph{Dirac series of $E_{8(-24)}$}, Int. Math. Res. Not. IMRN \textbf{2025}, rnaf058.

\bibitem{HC} Harish-Chandra, \emph{Harmonic analysis on real reductive groups I. The theory of the constant term}, J. Func. Anal. \textbf{19} (1975), 104--204.


\bibitem{HJ} S.~Helgason and K.~Johnson, \emph{The bounded spherical functions on symmetric spaces}, Adv. Math. \textbf{3} (1969), 586--593.


\bibitem{HP} J.-S.~Huang and P.~Pand\v zi\'c, \emph{Dirac
cohomology, unitary representations and a proof of a conjecture of
Vogan}, J. Amer. Math. Soc.  \textbf{15} (2002), 185--202.

\bibitem{HPV} J.-S.~Huang, P.~Pand\v zi\'c and  D.~Vogan,
\emph{On classifying unitary modules by their Dirac cohomology}, Sci. China Math. \textbf{60}  (2017), 1937--1962.

\bibitem{Kn} A. Knapp, \emph{Lie Groups Beyond an Introduction}, Birkh\"auser, 1996.

\bibitem{SV} S.~Salamanca-Riba and D.~Vogan, {\it On the classification of unitary representations of reductive Lie
groups}, Ann. of Math. \textbf{148} (1998), 1067--1133.

\bibitem{Vog86} D.~Vogan, {\it The unitary dual of ${\rm GL}(n)$ over an Archimedean field}, Invent. Math. \textbf{83} (1986), 449--505.

\bibitem{Vog97} D.~Vogan, \emph{Dirac operators and unitary
representations}, 3 talks at MIT Lie groups seminar, Fall 1997.
\end{thebibliography}
\end{document}